\documentclass{amsart}

\usepackage{amssymb}
\newcommand{\R}{\mathbb{R}}

\newcommand{\ra}{\rightarrow}
\newcommand{\ve}{\varepsilon}
 \newcommand{\vs}{\vspace*{1mm}}
 \newcommand{\ds}{\displaystyle}

\renewcommand{\int}{{\rm int\,}}
\newcommand{\sbs}{\subseteq}
\newcommand{\sgn}{{\rm sgn}}
\newcommand{\sft}{{\sf T}}

 \newtheorem{thm}{Theorem}[]
\newtheorem{theorem}{Theorem}[]
 
\newtheorem{corollary}[thm]{Corollary}
 
 \newtheorem{lemma}[thm]{Lemma}
 
 \theoremstyle{definition}
 
 \theoremstyle{remark}
 \newtheorem{remark}[thm]{Remark}
 
 \numberwithin{equation}{section}

\begin{document}

%-------------------------------------------------------------------------
% editorial commands: to be inserted by the editorial office
%
%\firstpage{1} \volume{228} \Copyrightyear{2004} \DOI{003-0001}
%
%
%\seriesextra{Just an add-on}
%\seriesextraline{This is the Concrete Title of this Book\br H.E. R and S.T.C. W, Eds.}
%
% for journals:
%
%\firstpage{1}
%\issuenumber{1}
%\Volumeandyear{1 (2004)}
%\Copyrightyear{2004}
%\DOI{003-xxxx-y}
%\Signet
%\commby{inhouse}
%\submitted{March 14, 2003}
%\received{March 16, 2000}
%\revised{June 1, 2000}
%\accepted{July 22, 2000}
%
%
%
%---------------------------------------------------------------------------
%Insert here the title, affiliations and abstract:
%

\title[Extrema and elementary symmetric functions]
 {A Note on Extrema of  Linear Combinations of Elementary Symmetric Functions}

%----------Author 1
\author{Alexander Kova\v{c}ec}

\address{%
Departamento de Matem\'atica \\
Faculdade de Ci\^encias e Tecnologia\\
Apartado 3008\\
EC Universidade\\
P 3001-454 Coimbra\\
Portugal}

\email{kovacec@mat.uc.pt}

%\thanks{This work was completed with the support of our
%\TeX-pert.}
%----------Author 2
\author{Salma Kuhlmann}
\address{Fachbereich Mathematik und Statistik \\
Universit\"at Konstanz \\
78457 Konstanz, Germany} \email{salma.kuhlmann@uni-konstanz.de}

%----------Author 3
\author{Cordian Riener}
\address{ Institut f\"ur Mathematik\\
Goethe Universit\"at \\
60325 Frankfurt \\
Germany}
\email{riener@math.uni-frankfurt.de}

%----------classification, keywords, date


\date{February 23, 2011}
%----------additions
%\dedicatory{To my boss}
%%% ----------------------------------------------------------------------

\begin{abstract}
This note provides a new approach to a result of Foregger \cite{for}
and related earlier results by Keilson \cite{keil} and Eberlein
\cite{eberl}. Using quite different techniques, we prove a more
general result from which the others follow easily. Finally, we
argue that the proof in \cite{for} is flawed.
\end{abstract}
\maketitle

\noindent Throughout this note write
$E_j(x)=E_j(x_1,...,x_n),$ $j=0,1,\ldots, n,$ for the $j$-th
elementary
 symmetric polynomial
in $n$ variables. Let $H=H_\gamma=\{x\in \R^n: E_1(x)=\gamma\}$ be a
hyperplane normal to the vector $1_n=(1,1,...,1)\in \R^n$ and $\phi$
be a real linear combination of the $E_j.$ Putting $D:= H \cap
[0,1]^n,$ Keilson \cite{keil} investigated the question of where the
function $\phi:D \ra \R$ assumes its extreme values. Eberlein
\cite{eberl} noted that the result of Keilson was actually already
known by Chebyshev in 1846 and Hoeffding in 1956. He assumed $0\leq
a_i\leq b_i\leq 1$ and investigated the same question on more
general domains $D'=H\cap ([a_1,b_1]\times \cdots \times
[a_n,b_n]).$ Foregger, interested in solving a problem by Pierce
\cite{pier}, also solved by Li \cite{ckli}, returned to the original
Keilson-Chebyshev question, and proposed to show in \cite[Theorem
2]{for} that if $\phi$ is non-constant on $D$ (i.e. has degree at
least two) and attains its extremum at an interior point of $D,$
then this point must be the symmetric point of $D,$ that is
$\frac{\gamma}{n} 1_n.$ At the end of his paper Foregger invites to
establish similar results to his for Eberlein's domains.

In this note, we establish below in Theorem \ref{theorem1}, Theorem
3 and Corollary 5 results from which the above mentioned are easily
derived. Our method, originally developed by Riener \cite{rien} to
provide an algebraic approach to Timofte's results \cite{tim} is
different and much simpler: it is based on the observation that the
elementary symmetric functions are closely related to roots of
univariate {\it hyperbolic polynomials}, i.e. polynomials with only
real roots. Given a function $f:D\sbs \R^n \ra \R,$  write
loc.extr$(f)$ for the set of local extrema of $f;$ that is for the
set of all points $u$ in $D$ admitting a neighborhood $N=N_u\sbs D$
such that  $\forall x\in N \;f(x)\geq f(u)$ or $\forall x\in N
\;f(x)\leq f(u).$
\begin{theorem} \label{theorem1}
If a  real linear combination  $\phi$ of elementary symmetric polynomials
has degree at least two, then it has as a function $\phi:H\ra \R$ on the
hyperplane $H=\{x: E_1(x)=\gamma\}$  at most one local extremum; more precisely
 $\mbox{\rm  loc.extr}(\phi)\sbs \{\frac{\gamma}{n}1_n\}.$
\end{theorem}
\noindent We will need an auxiliary lemma. \noindent
\begin{lemma}
Let $f\in \R[t]$ be hyperbolic of degree $n$ with
 only $s\in \{1,...,n-1\}$   {\em distinct} (real)
 roots all lying in some open interval $I.$ Then  there exists  $g\in \R[t]$ of degree $n-s$
 so that for all small $\ve>0$, the polynomials $f+ \ve g$ and $f- \ve g$  are hyperbolic
 and have both more than $s$ distinct roots all lying in $I.$
\end{lemma}
\begin{proof}Let $x_1, \ldots , x_s,$ be the distinct zeros of $f$  which by hypothesis
 all lie in $I.$ Without loss of generality  assume $I$ bounded and for notational
convenience assume
 $I=]x_0,x_{s+1}[$  and
  $x_0<x_1< \ldots <x_s<x_{s+1}.$
We have a factorization
       $$f=\prod_{i=1}^s (t-x_i)\cdot g_1=:p\cdot g_1, $$
where the roots of $g_1$ are all in  $\{x_1, \ldots, x_s\}$ and $g_1$ is of
degree $n-s\geq 1.$ Now
 for $x_i< t <x_{i+1},$  $i=0,\ldots,s,$  $\sgn\, p(t)=(-1)^{s-i}.$
Therefore there exist points $\xi_i\in ]x_i,x_{i+1}[,$ and $\ve>0,$ such that
                 if $s-i$ is odd, $p(\xi_i)<-2\ve,$ while if
$s-i$ is even, $p(\xi_i) >2\ve.$ Therefore the
polynomials $p\pm\ve$ have in
    $\xi_0,$
$\xi_1,\xi_2,...,\xi_{s-1},\xi_{s}$ alternating signs. So they have
in each of the intervals $]\xi_i,\xi_{i+1}[, $ $i=0,1,..., s-1,$ a
root. Therefore $p\pm\ve$ both have $s$ distinct
 real roots   and so are hyperbolic and   the roots all lie in $I$.    \quad
 Furthermore we see that
$p\pm \ve$ have none of their roots in the set $\{x_1, . . . ,
x_s\}$. Hence $(p \pm \ve) \cdot g_1 = f \pm \ve g_1$ are hyperbolic
and each have more than $s$ distinct roots all lying in $I$.
 \end{proof}
\begin{proof}[Proof of Theorem \ref{theorem1}] If  the set of  local extrema   of $\phi:H\ra \R$ is nonempty,
  choose a local extremum, say $a:=(a_1,\cdots,a_n)$,  so that
 the number $s:=|\{a_1,\ldots ,a_n\}|$ of distinct  components of $a$ is maximal  and   consider the
univariate polynomial \vs\centerline{ $\ds
f_a:=\prod_{i=1}^n(-a_i+t) =\sum_{i=0}^{n-2} (-1)^{n-i} E_{n-i}(a)
t^i  -\gamma t^{n-1}+ t^n.$ } Here we  used Vi\`ete's formula and to
facilitate reading  of the following we think of $f_a$ and similar
univariate polynomials as written from left to right with rising
powers and write $\phi$ as
$\ds\phi(x)=\sum_{i=0}^{n-2}c_{n-i}E_{n-i}(x)+ c_1 E_1(x)+c_0
E_0(x).$
 \noindent Assume  $s\in \{2,3,...,n-1\}.$ By lemma 2 there
exists a polynomial  $g=\sum_{i=0}^{n-s} g_{n-i} t^i$ of degree
$n-s,$ and an $\ve_0>0$ such that for all $\ve$ with
$0<|\ve|<\ve_0,$ the polynomial $f_a+\ve g$ is hyperbolic and has at
least $s+1$ distinct roots. We can represent for every $\ve$ the $n$
(real) roots of $f_a+\ve g$ as the entries of
$a(\ve)=(a_1(\ve),...,a_n(\ve)),$ and according to Kato
 \cite[p.109]{kato} think of $a(.)$ actually as a continuous
function. We then have a similar equation as above for $f_a+\ve g$
(in place of $f_a$) and $E_{n-i}(a(\ve))$ (in place of
$E_{n-i}(a)$). We find  $(-1)^{n-i}E_{n-i}(a(\ve))= (-1)^{n-i}
E_{n-i}(a)+  \ve g_{n-i},$ where for $i>n-s,$ we put $g_{n-i}=0.$ In
particular $E_1(a(\ve))=\gamma.$ By these formulae we find
\vs\centerline{$ \ds\phi(a(\ve))=\phi(a)+\ve
\sum_{i=0}^{n-s}(-1)^{n-i} c_{n-i}g_{n-i}.$} \noindent It now
follows for all  $\ve$ of small modulus and appropriate sign, that
$\phi(a(-\ve))\leq \phi(a)\leq \phi(a(\ve))$ with $a(\pm \ve)\in
O_a\sbs H,$ $O_a$ a neighborhood of $a.$
 So in every neighborhood of $a$ there are  points $a(-\ve)$
 at which $\phi$ is not larger than $\phi(a)$ and
 $a(\ve)$
 at which $\phi$ is not smaller than $\phi(a)$ and
 which have at least $s+1$ distinct coordinates.
This contradicts our choice of $a.$ Therefore $s\in \{1,n\},$ that
is, $f_a$ has either one root of multiplicity $n,$ or $n$ distinct
roots. Assume $s=n.$ Since degree$(\phi)\geq 2,$
$(c_2,...,c_n)\neq 0.$ Thus
 there exist reals $g_2,...,g_n$
so that $\sum_{i=0}^{n-2} c_{n-i}g_{n-i} \neq 0.$
Consider the polynomial $g=\sum_{i=0}^{n-2} g_{n-i} t^i$
and note that $f_a+\ve g$ will have for all $\ve$ of small modulus  $n$ roots.
With the reasoning above, we infer this time  strict inequalities
$\phi(a(-\ve))< \phi(a)< \phi(a(\ve)),$
  again arriving at a contradiction.
Hence $s=1$ and  $a=\frac{\gamma}{n}1_n,$ as we wished to show.
\end{proof}

\cite[Theorem 2]{for} follows at once from Theorem \ref{theorem1}:
If $\phi$ is non constant on $D$ then it is nonconstant on the
hyperplane $H$ of our theorem, hence evidently must have degree at
least two (that is some of the polynomials $E_2,...,E_n$ must occur
in $\phi$). If $\phi$ attains its maximum or minimum at an interior
point $P$ of $D$ then a small enough neighbourhood of $P$ as a point
of $D$ is also a neighbourhood of $P$ as a point in $H.$ So $P$ is a
local extremum as well for the function $\phi: D \ra R$ as for its
extension $\phi: H \ra R.$ Hence by our theorem $P$ must be the
symmetric point as claimed by Foregger.

A simple adaptation of the proof of Theorem \ref{theorem1}, which we
leave to the reader, yields the following generalization of Theorem
1 which will not be further used in this paper.
\begin{theorem} If the polynomial
 $\phi$ of theorem 1 has
 degree at least $k+1$ and is considered  as a function on the  real  variety
 $\{x: E_1(x)=\gamma_1, ..., E_k(x)=\gamma_k\},$
 then
each of its   local extrema
$a$ has at most $k$ distinct components, i.e.
$|\{a_1,\ldots, a_n\}|\leq k.$
\end{theorem}
\noindent We can now derive a corollary in the spirit of Eberlein
and Foregger which is  more complete. The following notation is
fixed through the rest of the paper. Let $a_i<b_i,$ $i=1,\ldots, n$
be real numbers. Denote by $e_i=(0,...,0,1,0,...0),$ the $i$-th
standard vector in $\R^n.$ Using the $(2n+2)\times n$ matrix

\centerline{$A=[-e_1^\sft,e_1^\sft, -e_2^\sft, e_2^\sft,\ldots,
-e_n^\sft,e_n^\sft, -1_n^\sft, 1_n^\sft]^\sft$ }

\noindent (whose rows are $-e_1, e_1, -e_2, ...$) and the column

\centerline{$d=(-a_1,b_1,-a_2,b_2,\ldots,-a_n,b_n,
-\gamma,\gamma)^\sft,$ }

\noindent the set
 $D'=H\cap ([a_1,b_1]\times \cdots \times [a_n,b_n])$ can be defined as $D'=\{x: Ax \leq d\}.$
Since $D'$ is  bounded it is a polytope in the sense of Schrijver \cite[p. 89]{sch}.
 For $K\sbs \{1,...,n\}$ we write $K^c=\{1,\ldots, n\}\setminus K$
for its complement. We begin by characterizing the faces of $D'.$
\begin{lemma}
A subset $F$ of $\R^n$ is  face of $D'$ if and only if it is
nonempty and there exists a subset $K$  in $\{1,2,\ldots, n\}$ and
$u_l\in \{a_l,b_l\}$ for all $l \in K^c,$  so that
\vs\centerline{$\ds F=\{\sum_{k\in K} x_ke_k +\sum_{l\in K^c}
u_le_l: \forall k\in K \, a_k\leq x_k
 \leq b_k   \mbox{ and } \sum_{k\in K} x_k+ \sum_{k\in K^c} u_k=\gamma \}.$}
\noindent If $K^c$ is a maximal set such that the inequality conditions at the right of the `:' are simultaneously
satisfiable without equalities, then the affine dimension of $F$ is $(|K|-1)^+.$
\end{lemma}
\begin{proof}
By \cite[p101]{sch}, a  set  $F\sbs \R^n$ is a face of $D'$ iff $F\neq \emptyset$ and
$F=\{x: Ax\leq d \mbox{ and } A'x=d'\},$
 where $A'$ and $d'$  are obtained from $A,d,$ respectively, by selecting the same row
indices. It may happen that there exists a further, say $i$-th  row $a_{i,*}$ of $A,$ such
that  $a_{i,*}x=d_i$  for all $x\in F.$ We think from now on of
 all such rows  as included in $A'$ and write $A^{=}$ for this submatrix of $A,$
and $d^{=}$ for the corresponding subvector of $d.$ Let $A^+$ be the
complementary submatrix of $A^{=}$ in $A.$ Then for each row index
$i$ of $A$ defining a row of $A^+$ there exists a vector $x^{(i)}\in
F$ such that $A^+x^{(i)}<d_i.$ The arithmetic mean of all such
vectors yields a vector  $x\in F$ so that $A^+x<d^+$ (meaning that
all component inequalities are strict). We note that  $F$ can be
written as $F=\{x: A^+ x\leq d \mbox{ and } A^{=}x=d^{=}\}.$ The set
of row indices of $A$ defining $A^{=}$ contains evidently $\{2n+1,
2n+2\}$ but contains   from each of the sets $\{2i-1,2i\},$
$i=1,...,n$ at most one index, since otherwise $F$ would be empty.
An index lying in $\{2i-1,2i\}$ fixes $x_i$ to be equal to some
$u_i\in \{a_i, b_i\}.$ Let $K^c\sbs \{1,...,n\}$ be the set of
indices  $i$  of variables so obtained, viewed as the
 complement of some set $K\sbs \{1,...,n\}.$  Then the components  $x_k,\; k\in K$ of $x\in F$
 may satisfy    $a_k<x_k < b_k$  and
$\sum_{k\in K} x_k=\gamma-\sum_{k\in K^c} u_k.$
It is now clear that $F$ can be written as above and the
affine dimension of $F$ (which is by definition the dimension of its affine hull)
equals $(|K|-1)^+.$
\end{proof}
\begin{corollary}
Assume that  as a function, $\phi:D'\ra \R$ is nonconstant on every
face of dimension one (edge) of the polytope $D'.$  Then every local
extremum $p$ of $\phi$
 is  point in a cartesian product
 $\{a_1, e, b_1\} \times \cdots \times \{a_n, e, b_n\}.$ Here $e$ has to be chosen such that
the sum of the components of $p$ is $\gamma.$
\end{corollary}
\begin{proof}
Let $p\in D'$ be a local extremum of $\phi.$ Since the faces of $D'$
form a lattice \cite[section 8.6]{sch}, there exists a unique face
$F$ of minimal dimension containing $p.$ Then $p\in\int F.$  We
think of $F$ as written in the above lemma, with the set $K^c$ there
chosen maximal. If  $\dim F=0,$ then  $p$ is a vertex,  $F=\{p\},$
and $|K|=0$ or $|K|=1.$ It then follows by the above
characterization of $F$ that $p=(p_1,\ldots, p_n)$ is a point so
that for
  at least $n-1$ indices $i,$  $p_i\in \{a_i,b_i\}.$
Consequently  $p$ lies  in one of the  cartesian
products admitted. Assume now $\dim F=k\geq 1.$ Then $F$ has a nonempty relative interior given via
$\gamma^*= \gamma-\sum_{k\in K^c} u_k$  by
\vs\centerline{$\ds \int F=\{\sum_{k\in K} x_ke_k +\sum_{l\in K^c} u_le_l: \forall k\in K \, a_k < x_k
 < b_k   \mbox{ and } \sum_{k\in K} x_k=\gamma^* \}.$}
\noindent Now $F$ has faces of $D'$ which have dimension 1 \cite[section 8.3]{sch} on which
by assumption $\phi$
 is not constant;
consequently $\phi$ is not constant on the affine hull of $F,$ and hence is a polynomial
of degree$\geq 2$ on it.
The local extremum $p\in \int F$ of $\phi$ is necessarily also local extremum of the restriction
 of $\phi$ to the affine hull of $F.$
Now this restriction is
$\phi'(x_{k_1},...,x_{k_{|K|}})=\phi'(\sum_{k\in K} x_ke_k)
:=\phi(\sum_{k\in K} x_ke_k +\sum_{l\in K^c} u_le_l).$ It is again a
certain real linear combination of elementary symmetric functions of
variables $x_k,$ $k\in K=\{k_1,...,k_{|K|}\}$: this  follows since a
similar fact holds for each of the elementary symmetric functions.
Applying our main result to the affine hull of  $F,$ we get that
 $p=\sum_{k\in K} \frac{\gamma^*}{|K|} e_k +\sum_{l\in K^c} u_le_l,$ as we wished to show.
\end{proof}
\begin{remark}
a. The mentioned authors state some of their results in somewhat
unprecise terms.
%Thus in \cite[p.269]{keil}
 %we find the statement that
%`$\phi$ assumes its maximum and minimum on the set $D$ at either its
%boundary or the symmetric point.' If one analyzes Keilson's proof
%one finds that he establishes the following: if $a\in D$ is a
%critical point of $\phi,$ then $a$ is the symmetric point or there
%is a line segment in $D$ on which $\phi$ is constant. It follows
%that if $a\in \mbox{int } D$ (the relative interior of $D$) is a
%{\it strict} local extremum  (i.e. a point so that for all points $
%x\in \mbox{int } D \setminus\{a\}$ lying in a neighbourhood of $a,$
%there holds $\phi(x)>\phi(a)$ or for all such $x$ there holds
%$\phi(x)<\phi(a)$), then $a$ is the symmetric point.  For Keilson's
%reasoning to be convincing his `extremum' on p. 270, line -12
%must be assumed to be strict.
 Eberlein \cite[Theorem 1 p. 312]{eberl} tells us that the
`minimum and the maximum of $\phi$ is assumed at least among the
points whose components which are not endpoints are all equal.
Moreover if the maximum and the minimum is attained only in the
interior of $D'$ then it is assumed uniquely at the point
$\frac{\gamma}{n}1_n$'.  After a  more precise formulation, the
corresponding proof is right but unfortunately does not give more.
In particular it does not exclude for nonconstant $\phi$ lines on
which $\phi$ has constant extreme value and which extend to the
boundary of $D'.$ Thus again there may be interior extrema other
than $\frac{\gamma}{n}1_n.$ This possibility is what Foregger
proposed to exclude for his domain $D.$ He attempts establishing his
result by induction over dimension via reasoning not easily
adaptable to Eberlein's more general domain $D'.$ We now discuss the
error in Foregger's proof, and point out another mistake in his
paper.

b. (i) Using exactly his notation, in \cite[p.384]{for}, Foregger
derives for $s<n,$ $c\in \R^{n-s}$ constant, and $y\in \R^s,$
 the first two lines in the following chain; the third line is a consequence of noting that $E_0(y)=1,$
 and $E_1(y)=\gamma^*$  by Foregger's (local) definition (of)
 $C_{\gamma^*}=\{y\in \R^s: \sum_{i=1}^s y_i= \gamma^*, y_i\in [0,1]\}\subset \R^s$ on p. 383:
\begin{align*}
  0&=\phi^*(y)\\
   &=\sum_{k=0}^n E_k(y)\sum_{r=k}^n c_r E_{r-k}(c) - \sum_{r=2}^n c_r E_r(0,c)\\
             &=\sum_{r=0}^n c_r E_r(c)+\gamma^*\sum_{r=1}^n c_r E_{r-1}(c) -
   \sum_{r=2}^n c_r E_r(0,c) + \sum_{k=2}^n E_k(y) \sum_{r=k}^n c_r E_{r-k}(c).
\end{align*}
Note that
 $y\in \R^s,$ so  the
 definition of the elementary symmetric functions requires $E_{s+1}(y)= \cdots =E_n(y)=0,$ a fact
not observed in \cite{for}. Granted that as functions $1,
E_2,E_3,..., E_s$ are (usually) linearly independent
 on $C_{\gamma^*}$ - to see this put $n-1$ variables equal to real variable $t,$ the other equal
to $\gamma^*-t$ and observe that  $E_j(t,...,t,\gamma^*-t)$ is a
polynomial of degree $j$ in $t$ - we  may infer the equations
$\ds\sum_{r=k}^n c_r E_{r-k}(c)=0, k=2,3,...,s.\>$\noindent The
problem in \cite{for} is that these equations are claimed for
$k=2,...,n,$ (and not only for $k=2,\ldots,s$) and then used in the
order $k=n,n-1,...,1,$ to derive that $c_n,c_{n-2},...,c_1,$ are 0.
Therefore the proof seems to be beyond repair.

(ii) Foregger in examining Eberlein's theorem, claims on  p. 385
that the function $\phi=\phi(x,y,z)=xyz-0.5( xy+xz+yz)$ assumes on
$C_{5/4}=\{(x,y,z):x+y+z=5/4\}\cap ([3/8,5/8]\times[3/8,5/8]\times
[1/8,3/8])$ in $p_0=(1/2,1/2,1/4)$ an interior maximum of value
$-0.1875=\phi(p_0).$ It is easily seen that this would contradict
his own (and our) main theorem. Indeed numerical experiments
indicate that $p_0$ is not an interior local extremum.

c. In \cite{waterh}, Waterhouse gives examples showing that the
symmetric point  in general needs not be  a local  extremum of
symmetric functions subject to symmetric conditions. We thank the
referee to have pointed out that $\phi(x)= E_2(x)-E_3(x)$ is an
example of a function that satisfies the hypothesis of our Theorem
\ref{theorem1} with $n=\gamma=3$ but has no extremum at all. Indeed
$\mathbb{R}\ni t\mapsto p(t)=(1+t,1+t,1-2t)$ describes a line in $H$
through the point $p(0)=(1,1,1)$ and $\phi(p(t))=2+2t^3.$  Examples
for larger $n$ can be similarly constructed.

\end{remark}
{\it Acknowledgment:} We thank Ant\'onio Leal Duarte for having
drawn our attention to Kato's book.

\bibliographystyle{gLMA}

\nocite{*}
 \end{document}